%% file: link-symplectic1.0.tex
\documentclass{amsart}
\include{headers}

\begin{document}
\title{Linked alternating forms and linked symplectic Grassmannians}
\author{Brian Osserman}
\author{Montserrat Teixidor i Bigas}
\begin{abstract} 
Motivated by applications to higher-rank Brill-Noether theory and the
Bertram-Feinberg-Mukai conjecture, we introduce the concepts of linked
alternating and linked symplectic forms on a chain of vector bundles,
and show that the linked symplectic Grassmannians parametrizing
chains of subbundles isotropic for a given linked symplectic form has
good dimensional behavior analogous to that of the classical symplectic
Grassmannian.
\end{abstract}

\maketitle

\section{Introduction}

Higher-rank Brill-Noether theory studies moduli spaces for pairs $(\sE,V)$ 
on a (smooth, projective) curve $C$, where $\sE$ is a vector bundle
of specified rank and degree, and $V$ is a vector space of global
sections of $\sE$, having specified dimension. The classical case of
line bundles is now well understood, with most proofs of the main theorems
relying on degeneration techniques. In the higher-rank case, there are a
number of partial results, with many of them using the generalization of 
the Eisenbud-Harris theory of limit linear series to higher-rank vector 
bundles given in \cite{te1}. However, even in the case of rank 
$2$, we have no comprehensive conjectures on the dimensions of the 
components of the moduli spaces, or when they are nonempty. See \cite{g-t1}
for a survey.

One phenomenon
observed by Bertram, Feinberg and Mukai in \cite{b-f2} and \cite{mu2}
is that loci of bundles of rank $2$ and canonical determinant always have
larger than the expected dimension, due to additional symmetries in this
case. This was generalized to other special determinants in \cite{os16}. 
In the case of canonical determinant, Bertram, Feinberg and Mukai 
conjectured that the moduli spaces in question were always
nonempty when their modified expected dimension was nonnegative, and
this conjecture remains open.
In order to use limit linear series techniques to prove existence results 
in this setting, it is necessary to understand how the symmetries arising
from special determinants interact with the expected dimension bounds of
the moduli spaces of generalized limit linear series. In the generalized
Eisenbud-Harris setting, this is not at all obvious, and our motivation
is to prove that we obtain the necessary modified expected dimension
bounds, using the alternative construction of limit linear series spaces 
presented in \cite{os8}. We emphasize that such results have immediate
implications: \cite{te5} proves rather strong existence results towards
the Bertram-Feinberg-Mukai conjecture assuming that the desired dimension
bounds hold for families of limit linear series on degenerations to 
chains of elliptic curves.

In the limit linear series construction of \cite{os8}, spaces of linked 
Grassmannians are introduced to serve as ambient moduli spaces. Given a
base scheme $S$, integers $r<d$, let $\sE_{\bullet}$ denote a chain of 
vector bundles $\sE_1,\dots,\sE_n$ on $S$ of rank $d$, together with 
homomorphisms 
$f_i:\sE_i \to \sE_{i+1}$ and $f^i:\sE_{i+1} \to \sE_i$ satisfying
certain natural conditions (recalled in Definition \ref{defn:basic}
below). Then the associated linked Grassmannian $LG(r,\sE_{\bullet})$
parametrizes tuples of subbundles $\sF_i \subseteq \sE_i$ of rank $r$,
which are all mapped into one another under the $f_i$ and $f^i$. These
schemes behave like flat degenerations of the classical Grassmannian 
$G(r,d)$, and indeed according to \cite{o-h1}, whenever the $f_i$ and $f^i$ 
are generically isomorphisms, the linked Grassmannian does in fact yield
a flat degeneration of $G(r,d)$. The basic idea of the construction of 
limit linear series spaces in \cite{os8} is that one replaces the 
Grassmannian used in the construction of linear series spaces on smooth
curves with a linked Grassmannian. On the other hand, one may express 
the modified expected dimension observed by Bertram, Feinberg and Mukai
by saying that one replaces the Grassmannian by a symplectic Grassmannian.
In order to combine the two, we thus wish to introduce a notation of
linked symplectic Grassmannian, and to prove that it has good dimension
behavior.

We first introduce a more general notion of a linked alternating
Grassmannian, based on a definition of linked alternating form. It turns
out that the key idea is to not only consider alternating forms on each
of the $\sE_i$, but also pairings between $\sE_i$ and $\sE_j$ for $i \neq j$,
satisfying certain natural compatibility conditions. We prove:

\begin{thm}\label{thm:main-alt} If 
$LAG(r,\sE_{\bullet},\left<,\right>_{\bullet}) 
\subseteq LG(r,\sE_{\bullet})$ is a linked alternating Grassmannian,
and $z \in LAG(r,\sE_{\bullet},\left<,\right>_{\bullet})$ is a smooth
point of $LG(r,\sE_{\bullet})$, then locally at $z$, we have that
$LAG(r,\sE_{\bullet},\left<,\right>_{\bullet})$ is cut out by $\binom{r}{2}$
equations inside $LG(r,\sE_{\bullet})$.
\end{thm}

We then define linked symplectic forms to be linked alternating forms 
satisfying a certain nondegeneracy condition, and thus define linked
symplectic Grassmannians to be the corresponding special case of linked
alternating Grassmannians. Via an analysis of tangent spaces, we prove:

\begin{thm}\label{thm:main-symp} If
$LSG(r,\sE_{\bullet},\left<,\right>_{\bullet}) 
\subseteq LG(r,\sE_{\bullet})$ is a linked symplectic Grassmannian,
and $S$ is regular,
and $z \in LSG(r,\sE_{\bullet},\left<,\right>_{\bullet})$ is a smooth
point of $LG(r,\sE_{\bullet})$, then $z$ is also smooth point of
$LSG(r,\sE_{\bullet},\left<,\right>_{\bullet})$, and the latter has 
codimension $\binom{r}{2}$ in $LG(r,\sE_{\bullet})$ at $z$.
\end{thm}

For the statements of the above theorems, note there is an explicit 
description of the smooth points of linked Grassmannians, recalled in
Theorem \ref{thm:lg} below.

The restriction to smooth points does not cause problems in proving 
existence results using limit linear series, as smoothing arguments based 
on the Eisenbud-Harris theory already require restriction to an open subset 
of the moduli space of limit linear series which is contained in the smooth
points of the ambient linked Grassmannian. The main limitation of our
result is then that the construction of \cite{os8} has thus far only
been carried out for reducible curves with two components. We are thus
able to conclude that limit linear series arguments are valid for pairs
with canonical determinant and degenerations to curves with two components;
see Theorem \ref{thm:canon-det-smoothing}. However,
in order to obtain the same result for arbitrary curves of compact
type, as needed for \cite{te5}, it is necessary to generalize the
construction of \cite{os8} to arbitrary curves of compact type.

\subsection*{Acknowledgements} 
We would like to thank Gavril Farkas for bringing to our attention
the problem of combining the dimension bounds for the canonical 
determinant and limit linear series constructions.

\section{Preliminaries}

We work throughout over a fixed base scheme $S$. We being by recalling
some ideas from \cite{os8} and \cite{os17}, rephrased in a more convenient
manner.

The basic definition is as follows:

\begin{defn}\label{defn:basic} Let $d,n$ be positive integers. 
Suppose that $\sE_1,\dots,\sE_n$ are vector bundles of rank $d$ on $S$
and we have homomorphisms 
$$f_i:\sE_i \to \sE_{i+1},\quad f^i:\sE_{i+1} \to \sE_i$$
for each $i=1,\dots,n-1$. Given $s \in \Gamma(S,\sO_S)$, we say that 
$\sE_{\bullet}=(\sE_i,f_i,f^i)_i$ is {\bf $s$-linked} if the 
following conditions are satisfied:
\begin{Ilist}
\itm For each $i=1,\dots,n$,
$$f_i \circ f^i = s \cdot \id, \text{ and } f^i \circ f_i = s \cdot \id.$$
\itm On the fibers of the $\sE_i$ at any point with $s=0$, we have that for
each $i=1,\dots,n-1$,
$$\ker f^i = \im f_i, \text{ and } \ker f_i = \im f^i.$$
\itm On the fibers of the $\sE_i$ at any point with $s=0$, we have that for
each $i=1,\dots,n-2$,
$$\im f_i \cap \ker f_{i+1}=(0),\text{ and }\im f^{i+1} \cap \ker f^i = (0).$$
\end{Ilist}
If $\sE_{\bullet}$ satisfies conditions (I) and (III), we say it is
{\bf weakly $s$-linked}.
\end{defn}

The $s$-linkage condition is precisely that required in \cite{os8} for the 
ambient bundles for a linked Grassmannian, which we also recall below.
The following notation will be convenient:

\begin{notn} In the situation of Definition \ref{defn:basic}, with $i<j$ 
we write
$$f_{i,j}=f_{j-1} \circ \cdots \circ f_i$$
and
$$f^{j,i}=f^i \circ \cdots \circ f^{j-1}.$$
We also write $f_{i,i}=\id, f^{i,i}=\id$.
\end{notn}

We have the following basic structure for $s$-linked bundles:

\begin{lem}\label{lem:structure} Suppose that $\sE_{\bullet}$
is $s$-linked. Let $r_i=\rk f_i$ for $i=1,\dots,n-1$, and by convention
set $r_0=0$, $r_n=d$.
Then locally on $S$, for $i=1,\dots,n$ there exist subbundles 
$\sW_i \subseteq \sE_i$ of rank $r_i-r_{i-1}$ such that:
\begin{ilist}
\itm For $i=2,\dots,n-1$ we have that 
$$\sW_i \cap \spn(\ker f_i, \ker f^{i-1})=(0),$$ 
and similarly $\sW_1 \cap \ker f_1=(0), \sW_n \cap \ker f^{n-1}=(0)$.
\itm For all $j<i$, the restriction of $f_{j,i}$ to $\sW_j$ is an isomorphism
onto a subbundle of $\sE_i$, and for $j>i$ the restriction of $f^{j,i}$ to
$\sW_j$ is an isomorphism onto a subbundle of $\sE_i$.
\itm The natural map
$$\left(\bigoplus_{j=1}^{i} f_{j,i}(\sW_j)\right) 
\oplus \left(\bigoplus_{j=i+1}^{n} f^{j,i}(\sW_j)\right) \to \sE_i$$
is an isomorphism for each $i$. 
\end{ilist}
\end{lem}

The proof is similar to that of Lemma 2.5 of \cite{os17} (see also the
proof of Lemma A.12 (ii) of \cite{os8}), and is omitted.

We next discuss moduli of linked subbundles.

\begin{defn} Given $\sE_{\bullet}$ weakly $s$-linked of rank $d$,
and $r<d$, a {\bf linked subbundle} $\sF_{\bullet} \subseteq \sE_{\bullet}$
of rank $r$ consists of a subbundle $\sF_i \subseteq \sE_i$ for each $i$
such that $f_i \sF_i \subseteq \sF_{i+1}$, and $f^i \sF_{i+1} \subseteq
\sF_i$ for $i=1,\dots,n-1$. 
\end{defn}

Note that a linked subbundle automatically inherits the structure of a
weakly $s$-linked subbundle, but it is not necessarily the case that
a linked subbundle of an $s$-linked bundle is $s$-linked.

\begin{defn} Suppose $\sE_{\bullet}$ is $s$-linked of rank $d$, and we
are given $r<d$. Then the {\bf linked Grassmannian} 
$LG(r,\sE_{\bullet})$ is the scheme
representing the functor of linked subbundles of $\sE_{\bullet}$
of rank $r$.
\end{defn}

It is easy to see that $LG(r,\sE_{\bullet})$ exists, and is in fact a 
projective scheme over $S$, as it is cut out as a closed subscheme
of the product of classical Grassmannians
$$G(r,\sE_1) \times_S \cdots \times_S G(r,\sE_n).$$
Also note that for fibers with $s \neq 0$, condition (I) of $s$-linkage
implies that all the $f_i$ and $f^i$ are isomorphisms, so any subbundle
$\sF_i$ uniquely determines the others. Thus, the corresponding fiber of 
$LG(r,\sE_{\bullet})$ is isomorphic to the classical Grassmannian
$G(r,d)$. The interesting question is thus what happens at points with
$s=0$.

An important definition is:

\begin{defn} Given $\sE_{\bullet}$ $s$-linked on $S$, a morphism $T \to S$,
and a linked subbundle $\sF_{\bullet} \subseteq \sE_{\bullet}|_T$ of rank
$r$, we say that $\sF_{\bullet}$ is an {\bf exact point} of 
$LG(r,\sE_{\bullet})$ if
on the fibers of the $\sF_i$ at any point of $T$ with $s=0$, we have that for
each $i=1,\dots,n-1$,
$$\ker f^i = \im f_i, \text{ and } \ker f_i = \im f^i.$$
\end{defn}

Equivalently, $\sF_{\bullet}$ is an exact point if it is $s$-linked.
It is not hard to see that the exact points form an open subscheme of
$LG(r,\sE_{\bullet})$. The main results 
of \cite{os8} on linked Grassmannians, illustrating the value
of the $s$-linkage condition, are then the following:

\begin{thm}\label{thm:lg} If $\sE_{\bullet}$ is $s$-linked of rank $d$,
and we are given $r<d$, the exact points of $LG(r,\sE_{\bullet})$ are 
precisely the smooth points of $LG(r,\sE_{\bullet})$ over $S$. They have
relative dimension $r(d-r)$, and are dense in every fiber.
\end{thm}

These are Lemma A.12, Proposition A.13, and Lemma A.14 of \cite{os8}.

We thus see that the linked Grassmannian gives degenerations of the
classical Grassmannian (in fact, according to the main result of
\cite{o-h1} these are flat and Cohen-Macaulay degenerations, but this
will not be important for us). The following result, which is contained
in the proof of Lemma A.14 of \cite{os8}, will also be important:

\begin{lem}\label{lem:exact-tangent} In the case that $S$ is a point, let 
$\sF_{\bullet} \subseteq \sE_{\bullet}$ be an exact point of 
$LG(r,\sE_{\bullet})$, and let $(\sW_i \subseteq \sF_i)_i$ be as in Lemma 
\ref{lem:structure}. Then the tangent space to $LG(r,\sE_{\bullet})$ at
the point corresponding to $\sF_{\bullet}$ is canonically identified with
$\bigoplus_i \Hom(\sW_i,\sE_i/\sW_i)$.
\end{lem}

\section{Linked alternating forms}

We now introduce the definitions of linked bilinear form and linked
alternating form which will be central to our analysis.

\begin{defn}\label{defn:link-bilin} Given a weakly $s$-linked
$\sE_{\bullet}=(\sE_i,f_i,f^i)_i$, and $m \in \frac{1}{2} \ZZ$ 
between $1$ and $n$, a {\bf linked bilinear form} of {\bf index} $m$ on
$\sE_{\bullet}$ is a collection of bilinear pairings for each $i,j$
$$\left<,\right>_{i,j}:\sE_i \times \sE_j \to \sO_S$$
satisfying the following compatibility conditions: 
for all suitable $i,j$, we have
$$\left<,\right>_{i,j} \circ (f_{i-1} \times \id) = 
s^{\epsilon_{i,j}} \left<,\right>_{i-1,j},$$
$$\left<,\right>_{i,j} \circ (\id \times f_{j-1}) = 
s^{\epsilon_{i,j}} \left<,\right>_{i,j-1},$$
$$\left<,\right>_{i,j} \circ (f^i \times \id) = 
s^{\epsilon^{i,j}} \left<,\right>_{i+1,j}, \text{ and}$$
$$\left<,\right>_{i,j} \circ (\id \times f^j) = 
s^{\epsilon^{i,j}} \left<,\right>_{i,j+1},$$
where 
$\epsilon_{i,j}=\begin{cases}1 : i+j>2m \\ 0: \text{otherwise,}\end{cases}$
and 
$\epsilon^{i,j}=\begin{cases}1 : i+j<2m \\ 0: \text{otherwise.}\end{cases}$
\end{defn}

Note that the data determining a linked bilinear form can be described 
equivalently as a bilinear form on $\bigoplus_i \sE_i$, but it is harder
to describe the compatibility conditions in this context. 
We also observe that if
$m$ is an integer, then each fixed-index form $\left<,\right>_{i,i}$ is 
induced from $\left<,\right>_{m,m}$ by setting
$$\left<,\right>_{i,i}=\left<,\right>_{m,m}\circ (f_{i,m} \times f_{i,m})$$
for $i<m$, and similarly with $f^{i,m}$ for $i>m$. For further discussion of
the motivation for and consequences of the compatibility conditions, see 
Remark \ref{rem:compat} below. 

The following lemma checks that our compatibility conditions are internally
consistent, and will be useful later. Its proof is trivial from the
definitions.

\begin{lem}\label{lem:compat-consist} The compatibility conditions imposed
in Definition \ref{defn:link-bilin} satisfy the following 
consistencies:
\begin{ilist} 
\itm For all suitable $i,j$, they impose that
$$\left<,\right>_{i,j} \circ (f_{i-1} \times \id) \circ (f^{i-1} \times \id)
=s \left<,\right>_{i,j},$$
$$\left<,\right>_{i,j} \circ (f^i \times \id) \circ (f_i \times \id)
=s \left<,\right>_{i,j},$$
$$\left<,\right>_{i,j} \circ (\id \times f_{j-1}) \circ (\id \times f^{j-1})
=s \left<,\right>_{i,j},\text{ and}$$
$$\left<,\right>_{i,j} \circ (\id \times f^{j}) \circ (\id \times f_{j})
=s \left<,\right>_{i,j}.$$
Formally, 
$$\epsilon_{i,j}+\epsilon^{i-1,j}=1,$$
$$\epsilon_{i+1,j}+\epsilon^{i,j}=1,$$
$$\epsilon_{i,j}+\epsilon^{i,j-1}=1,\text{ and}$$
$$\epsilon_{i,j+1}+\epsilon^{i,j}=1.$$
\itm For all suitable $i,j$, they impose that
$$\left<,\right>_{i,j} \circ (f_{i-1} \times \id) \circ (\id \times f_{j-1})
=\left<,\right>_{i,j} \circ (\id \times f_{j-1}) \circ (f_{i-1} \times \id),$$
$$\left<,\right>_{i,j} \circ (f_{i-1} \times \id) \circ (\id \times f^{j})
=\left<,\right>_{i,j} \circ (\id \times f^{j}) \circ (f_{i-1} \times \id),$$
$$\left<,\right>_{i,j} \circ (f^{i} \times \id) \circ (\id \times f_{j-1})
=\left<,\right>_{i,j} \circ (\id \times f_{j-1}) \circ (f^{i} \times \id),
\text{ and}$$
$$\left<,\right>_{i,j} \circ (f^{i} \times \id) \circ (\id \times f^{j})
=\left<,\right>_{i,j} \circ (\id \times f^{j}) \circ (f^{i} \times \id).$$
Formally,
$$\epsilon_{i,j}+\epsilon_{i-1,j}=\epsilon_{i,j}+\epsilon_{i,j-1},$$
$$\epsilon_{i,j}+\epsilon^{i-1,j}=\epsilon^{i,j}+\epsilon_{i,j+1},$$
$$\epsilon^{i,j}+\epsilon_{i+1,j}=\epsilon_{i,j}+\epsilon^{i,j-1},\text{ and}$$
$$\epsilon^{i,j}+\epsilon^{i+1,j}=\epsilon^{i,j}+\epsilon^{i,j+1}.$$
\end{ilist}
\end{lem}

The alternating condition is then imposed as follows.

\begin{defn}\label{defn:link-alt} In the notation of Definition 
\ref{defn:link-bilin}, a linked bilinear form is a {\bf linked alternating
form} if $\left<,\right>_{i,i}$ is an alternating form on $\sE_i$ for all 
$i$, and $\left<,\right>_{i,j}=-\left<,\right>_{j,i} \circ \sw_{i,j}$ for all
$i \neq j$, where $\sw_{i,j}:\sE_i \times \sE_j \to \sE_j \times \sE_i$
is the canonical switching map.
\end{defn}

This definition is equivalent to requiring that the induced form on
$\bigoplus_i \sE_i$ be alternating.

Observe that being $s$-linked or weakly $s$-linked is preserved by
base change. It thus makes sense to define moduli functors of linked
bilinear forms and linked alternating forms (and indeed, one can do
this without any linkage conditions, if $s$ is given). Moreover, it
is clear that these functors have natural module structures, and
are representable. Our first result is that for $s$-linked bundles, the 
moduli of linked bilinear forms and of linked alternating forms behave 
just like their classical counterparts. 

\begin{prop}\label{prop:form-structure} Suppose 
$\sE_{\bullet}=(\sE_i,f_i,f^i)_i$ is
$s$-linked, and $m \in \frac{1}{2}\ZZ$ is between $1$ and $n$. Then
the moduli scheme of linked bilinear forms on $\sE_{\bullet}$ of index 
$m$ is a vector bundle on $S$ of rank $d^2$, and 
the moduli scheme of linked alternating forms on $\sE_{\bullet}$ of 
index $m$ is a vector bundle on $S$ of rank $\binom{d}{2}$.
\end{prop}

\begin{proof} First, choose subbundles $\sF_i \subseteq \sE_i$ as provided
by Lemma \ref{lem:structure}. Clearly, a linked bilinear form on
$\sE_{\bullet}$ induces by restriction a collection of bilinear pairings
$$\left<,\right>_{i,j}':\sF_i \times \sF_j \to \sO_S,$$
or equivalently a bilinear form on $\bigoplus_i \sF_i$,
and our claim is that this restriction map induces an isomorphism of 
functors from linked bilinear forms to bilinear forms on $\bigoplus_i \sF_i$. 
Because $\sum_i \rk \sF_i = d$, the claim yields the first statement
of the proposition. 

To prove the claim, suppose we have a collection of $\left<,\right>_{i,j}'$
as above; we aim to construct an inverse to the restriction map. Because 
$$\sE_i \cong \left(\bigoplus_{j=1}^{i} f_{j,i}(\sF_j)\right) 
\oplus \left(\bigoplus_{j=i+1}^{n} f^{j,i}(\sF_j)\right),$$
in order to define $\left<v_1,v_2\right>_{i,j}$, it is enough to do so for
$v_1$ in either $f_{\ell,i} \sF_{\ell}$ with $\ell \leq i$, or
$f^{\ell,i} \sF_{\ell}$ with $\ell >i$, and $v_2$ in either 
$f_{\ell',j} \sF_{\ell'}$ with $\ell' \leq j$, or
$f^{\ell',j} \sF_{\ell'}$ with $\ell' >j$. Starting from the necessity
of having $\left<,\right>_{i,j}=\left<,\right>_{i,j}'$ on
$\sF_i \times \sF_j$, we then see that inductive application of the 
compatibility conditions
of Definition \ref{defn:link-bilin} uniquely determine $\left<,\right>_{i,j}$.
We have to check that the resulting $\left<,\right>_{i,j}$ is well defined,
and satisfies all the compatibility conditions. This is straightforward to
verify case by case, using Lemma \ref{lem:compat-consist}.
The claim then follows, as the preceding construction is visibly inverse 
to the restriction map.

To obtain the second statement of the proposition, it is enough to observe
that under the isomorphism of functors constructed above, a linked
bilinear form is alternating if and only if the induced form on
$\bigoplus_i \sF_i$ is alternating. Indeed, this follows from the
symmetry of the compatibility conditions together with Lemma 
\ref{lem:compat-consist}.
\end{proof}

Proposition \ref{prop:form-structure} has immediate consequences for loci of 
isotropy. The relevant definitions are as follows.

\begin{defn}\label{defn:link-isotropy} If $\sE_{\bullet}$ is weakly
$s$-linked, with a linked bilinear form $\left<,\right>_{\bullet}$, 
we say that $\sE_{\bullet}$ is {\bf isotropic} for $\left<,\right>_{\bullet}$
if for all $i,j$, we have that $\left<,\right>_{i,j}$ vanishes uniformly.
\end{defn}

\begin{defn}\label{defn:link-isotropy-locus} If $\sE_{\bullet}$ is 
weakly $s$-linked, with a linked bilinear form $\left<,\right>_{\bullet}$, 
the {\bf locus of isotropy} of $\left<,\right>_{\bullet}$ on $\sE_{\bullet}$
is the closed subscheme of $S$ representing the functor of morphisms
$T \to S$ such that $\sE_{\bullet}$ is isotropic for 
$\left<,\right>_{\bullet}$ after restriction to $T$.
\end{defn}

The fact
that the locus of isotropy is represented by a closed subscheme is clear,
as $\sE_{\bullet}$ together with $\left<,\right>_{\bullet}$
induces a morphism from $S$ to the moduli scheme of linked bilinear forms
on $\sE_{\bullet}$, and the locus of isotropy is the preimage under
this morphism of the zero form.

Proposition \ref{prop:form-structure} thus implies:

\begin{cor}\label{cor:iso-codim} Suppose $(\sE_i,f_i,f^i)_i$ is
$s$-linked, and $(\left<,\right>_{i,j})_{i,j}$ is a linked bilinear
(respectively, linked alternating) form on $S$. Then the locus of
$S$ on which $(\sE_i)_i$ is isotropic is locally cut out by
$d^2$ (respectively, $\binom{d}{2}$) equations, and thus if $S$
is locally Noetherian, every component of this locus has codimension
at most $d^2$ (respectively, $\binom{d}{2}$) in $S$.
\end{cor}

\begin{rem}\label{rem:compat} We conclude with a discussion of the
motivation for Definition \ref{defn:link-bilin}. The idea, at least in the
case that $m \in \ZZ$, is that all of the forms are induced from a single
form $\left<,\right>_{m,m}$, which in our ultimate application will be
nondegenerate. In this situation, we cannot avoid having 
$\left<,\right>_{i,i}$ be degenerate on $\ker f_i$ for $i<m$ and $\ker f^i$
for $i>m$, and examples show that if we only consider the forms 
$\left<,\right>_{i,i}$, we will not obtain the behavior we want.
Because the $\left<,\right>_{i,i}$ are not uniformly zero, there is no way
to modify them to make them nondegenerate.

However, if suppose that $S=\Spec A$, with
$A$ a DVR, and $s \in A$ a uniformizer, and if we have a nondegenerate form
$\left<,\right>_{m,m}$, then the maps $f_{\bullet}$ and $f^{\bullet}$
induce not only forms $\left<,\right>_{i,i}$, but also pairings 
$\left<,\right>_{i,j}$ for all $i,j$. In the cases of interest to us however,
on the special fiber we will have $\im f_{m-1}$ orthogonal to $\im f^m$ in 
$\sE_m$, so if we simply take the induced pairings, we will have 
$\left<,\right>_{i,j}=0$ uniformly on the special fiber if $i<m$ and $j>m$,
or vice versa. But this means that considered over all of $S$, the forms
$\left<,\right>_{i,j}$ are multiples of $s$, and we can factor out powers
of $s$ (of exponent equal to $\min(|m-i|, |m-j|)$) so that the form does not
vanish uniformly on the special fiber. In the cases of interest to us, we 
will actually obtain nondegenerate forms this way when $i+j=2m$.

For an example of the importance of this additional nondegeneracy,
see Example \ref{ex:need-pairings}.
\end{rem}

\section{Linked symplectic forms}

Our next task is to give a suitable notion of nondegeneracy for linked
alternating forms, which we will use to define linked symplectic
Grassmannians. 

\begin{defn}\label{defn:link-symplectic} Suppose that 
$\sE_{\bullet}=(\sE_i,f_i,f^i)_i$ is
weakly $s$-linked, and $\left<,\right>_{\bullet}$ is a linked alternating
form on $\sE_{\bullet}$. We say that $\left<,\right>_{\bullet}$ is a 
{\bf linked symplectic form} if the following conditions are satisfied:
\begin{Ilist}
\itm for all $i,j$ between $1$ and $n$ with $i+j=2m$, we have
$\left<,\right>_{i,j}$ nondegenerate.
\itm if $2m<n+1$, then on all fibers where $s=0$, and for all $i$ with 
$2m-1<i \leq n$, the degeneracy
of $\left<,\right>_{i,1}$ is equal to $\ker f^{i-1}$.
\itm if $2m>n+1$, then on all fibers where $s=0$, and for all $i$ with 
$1\leq i < 2m-n$, the degeneracy
of $\left<,\right>_{i,n}$ is equal to $\ker f_{i}$.
\end{Ilist}
\end{defn}

Note that in conditions (II) and (III), the compatibility conditions of
Definition \ref{defn:link-bilin} imply that the degeneracy is at least
the specified subspaces, so all of the conditions are nondegeneracy
conditions, and we obtain an open subset of all linked alternating forms.

The following construction will be used to analyze the tangent space to the
linked symplectic Grassmannian.

\begin{defn}\label{defn:tangent-form} Suppose that 
$\sE_{\bullet}=(\sE_i,f_i,f^i)_i$ is $s$-linked, and 
$\left<,\right>_{\bullet}$ is a linked alternating form on $\sE_{\bullet}$. 
Let $\sF_{\bullet} \subseteq \sE_{\bullet}$ be an exact linked subbundle,
and suppose that $\sF_{\bullet}$ is isotropic for (the restriction of)
$\left<,\right>_{\bullet}$. Finally, let $(\sW_i \subseteq \sF_i)_i$ be as 
in Lemma \ref{lem:structure}. Given
a tuple of homomorphisms $(\vp_i:\sW_i \to \sE_i/\sW_i)_{1=1,\dots,n}$ 
define the
associated linked alternating form $\left<,\right>^{\vp_{\bullet}}_{\bullet}$
on $\sF_{\bullet}$ by applying the following formula on the $\sW_i$:
$$\left<,\right>^{\vp_{\bullet}}_{i,j} =
\left<,\right>_{i,j} \circ (\vp_i \times \id) + 
\left<,\right>_{i,j} \circ (\id \times \vp_j).$$ 
\end{defn}

Note that this is well-defined because $\sF_{\bullet}$ is assumed to be 
isotropic. Also, recall that by Proposition \ref{prop:form-structure},  
the pairings on the $\sW_i$ defined above uniquely determine a linked
alternating form $\left<,\right>^{\vp_{\bullet}}_{\bullet}$ on 
$\sF_{\bullet}$. We then have the following consequence of the symplectic
condition:

\begin{lem}\label{lem:symplectic-transverse} In the situation of
Definition \ref{defn:tangent-form}, suppose further that
$\left<,\right>_{\bullet}$ is a linked symplectic form, and that $S$ is
a point. Then the 
map from $\bigoplus_{i=1}^n \Hom(\sW_i,\sE_i/\sW_i)$ to the space of
linked alternating forms on $\sF_{\bullet}$ is surjective.
\end{lem}

For the proof of Lemma \ref{lem:symplectic-transverse}, the following
lemma is helpful. The proof is an immediate consequence of the compatibility
conditions of Definition \ref{defn:link-alt}.

\begin{lem}\label{lem:form-induced} Let $\sE_{\bullet}$ be $s$-linked,
with $\left<,\right>_{\bullet}$ a linked bilinear form of index $m$
on $\sE_{\bullet}$. Then:
\begin{ilist}
\itm given $i,j$ with $i+j>2m$, and any $\ell$ between $1$ and $n$ with 
$2m-i \leq \ell <j$, we have
$$\left<,\right>_{i,j}=\left<,\right>_{i,\ell} \circ (\id \times f^{j,\ell});$$
\itm given $i,j$ with $i+j<2m$, and any $\ell$ between $1$ and $n$ with 
$j < \ell \leq 2m-i$, we have
$$\left<,\right>_{i,j}=\left<,\right>_{i,\ell} \circ (\id \times f_{j,\ell});$$
\itm given $i,j$ with $i+j>2m$, and any $\ell$ between $1$ and $n$ with 
$2m-j \leq \ell <i$, we have
$$\left<,\right>_{i,j}=\left<,\right>_{\ell,j} \circ (f^{i,\ell} \times \id);$$
\itm given $i,j$ with $i+j<2m$, and any $\ell$ between $1$ and $n$ with 
$i < \ell \leq 2m-j$, we have
$$\left<,\right>_{i,j}=\left<,\right>_{\ell,j} \circ (f_{i,\ell} \times \id).$$
\end{ilist}
\end{lem}

Note in particular that if, for instance, $2m-i$ is between $1$ and $n$,
then $\left<,\right>_{i,j}$ is induced from $\left<,\right>_{i,2m-i}$.

We will also use the following easy lemma from linear algebra:

\begin{lem}\label{lem:subsetneq} Suppose $\left<,\right>:V \times W \to k$
is a non-degenerate bilinear pairing of $k$-vector spaces, and we have 
subspaces $W_1 \subsetneq W_2 \subseteq W$ and $V' \subseteq V$. If 
$(V')^{\perp} \cap W_2 \subseteq W_1$ in $W$, then 
$V' \cap W_2^{\perp} \subsetneq V' \cap W_1^{\perp}$ in $V$.
\end{lem}

\begin{proof}[Proof of Lemma \ref{lem:symplectic-transverse}] If we choose
bases $v^i_j$ for each $\sW_i$, it is clearly enough to prove that for all
$i,j,p,q$, unless $i=j$ and $p=q$ there exists a choice
of $\vp_{\bullet}$ such that the pairing
$\left<v^{i'}_{p'},v^{j'}_{q'}\right>^{\vp_{\bullet}}_{i',j'}$ is zero
for all $i',j',p',q'$ except $i'=i,j'=j,p'=p,q'=q$, or $i'=j,j'=i,p'=q,q'=p$,
and in these last two cases, the pairing is nonzero. Given $i,j,p,q$, 
first suppose $|m-i| \leq |m-j|$. Then set $\vp_{i'}=0$ for all $i' \neq i$,
and set $\vp_{i}(v^{i}_{p'})=0$ for all $p' \neq p$. We then wish to show
that there exists a choice of $\vp_i(v^i_p) \in \sE_i$ such that
$\left<\vp_i(v^i_p),v^j_q\right>_{i,j} \neq 0$, but 
$\left<\vp_i(v^i_p),v^{j'}_{q'}\right>_{i,j'}=0$ for all other choices of
$j',q'$. Equivalently, if we denote by $\hat{\sW}_{j} \subseteq \sW_j$ the 
span of the $v^j_{q'}$ for $q' \neq q$, we want 
$$\vp_i(v^i_p) 
\in (\hat{\sW}_j)^{\perp} \cap \left(\cap_{j'\neq j} \sW_{j'}^{\perp}\right),$$
but
$$\vp_i(v^i_p) \not\in \cap_{j'=1}^n \sW_{j'}^{\perp}.$$
Here each orthogonal space should be taken with respect to the appropriate 
pairing.

Now, if we have $1 \leq 2m-i \leq n$, then according to Lemma 
\ref{lem:form-induced}, the above conditions are equivalent to having
\begin{align*}\vp_i(v^i_p) &
\in (f \hat{\sW}_j)^{\perp} \cap 
\left(\bigcap_{j' \neq j, j' \leq 2m-i} (f_{j',2m-i} \sW_{j'})^{\perp} \right)
\cap
\left(\bigcap_{j' \neq j, j' > 2m-i} (f^{j',2m-i} \sW_{j'})^{\perp}\right) \\
& = \left(f(\hat{\sW}_j) \oplus 
\left(\bigoplus_{j' \neq j, j' \leq 2m-i} f_{j',2m-i} \sW_{j'}\right) \oplus
\left(\bigoplus_{j' \neq j, j' > 2m-i} f^{j',2m-i} \sW_{j'}\right)
\right)^{\perp},
\end{align*}
where $f=f_{j,2m-i}$ or $f=f^{j,2m-i}$ as appropriate, but
\begin{align*} \vp_i(v^i_p) & \not\in 
\left(\bigcap_{j'=1}^{2m-i} (f_{j',2m-i} \sW_{j'})^{\perp} \right) 
\cap \left(\bigcap_{j'=2m-i+1}^n (f^{j',2m-i} \sW_{j'})^{\perp} \right) \\
& = \left(\left(\bigoplus_{j'=1}^{2m-i} f_{j',2m-i} \sW_{j'}\right) 
\oplus \left(\bigoplus_{j'=2m-i+1}^n f^{j',2m-i} \sW_{j'}\right)
\right)^{\perp}.
\end{align*}
The sums are direct
sums because of Lemma \ref{lem:structure}, and now all the orthogonal 
complements are relative to $\left<,\right>_{i,2m-i}$. Again by Lemma
\ref{lem:structure}, the two sums
give distinct subspaces of $\sF_{2m-i}$, so by the nondegeneracy of
$\left<,\right>_{i,2m-i}$ imposed in the definition of a linked symplectic
form, we conclude that a $\vp_i(v^i_p)$ satisfying the desired conditions
exists in this case. 

On the other hand, if $2m-i < 1$, we can still apply Lemma 
\ref{lem:form-induced} to conclude that what we want is equivalent to
$$ \vp_i(v^i_p) \in 
\left(f^{j,1} (\hat{\sW}_j) \oplus 
\left(\bigoplus_{j' \neq j} f_{j',1} \sW_{j'}\right)\right)^{\perp},$$
but
$$ \vp_i(v^i_p) \not \in
\left(\bigoplus_{j'=1}^{n} f_{j',1} \sW_{j'}\right)^{\perp},$$
where the orthogonal complements are relative to $\left<,\right>_{i,1}$.
Now, $2m-1<i$, so applying Lemma \ref{lem:form-induced} again, what we
want is equivalent to
$$f^{i,2m-1} \vp_i(v^i_p) \in 
\left(f^{j,1} (\hat{\sW}_j) \oplus 
\left(\bigoplus_{j' \neq j} f_{j',1} \sW_{j'}\right)\right)^{\perp},$$
but
$$\vp_i(v^i_p) \not \in
\left(\bigoplus_{j'=1}^{n} f_{j',1} \sW_{j'}\right)^{\perp},$$
where now the orthogonal complements are relative to 
$\left<,\right>_{2m-1,1}$. Since this form is by the symplectic condition
nondegenerate, we have that there exist vectors in $\sE_{2m-1}$ with the
desired properties, and it is enough to show that we may further assume
they lie in $f^{i,2m-1}(\sE_i)$. To show this, by Lemma \ref{lem:subsetneq}
it is enough to show
\begin{equation}\label{eq:subset} (f^{i,2m-1} (\sE_i))^{\perp} \cap 
\left(\bigoplus_{j'=1}^{n} f_{j',1} \sW_{j'}\right) \subseteq 
f^{j,1} (\hat{\sW}_j) \oplus 
\left(\bigoplus_{j' \neq j} f_{j',1} \sW_{j'}\right).\end{equation}
We then observe that the
hypothesis that the degeneracy of $\left<,\right>_{i,1}$ on $\sE_i$ is
equal to $\ker f^i=\ker f^{i,2m-1}=\ker f^{i,1}$ implies that
$$\sE_{1}=(f^{i,2m-1}(\sE_i))^{\perp}\oplus f^{i,1}(\sE_i),$$
since the subspaces are of complementary dimension and have trivial 
intersection. Finally, we use the hypothesis that $|m-i| \leq |m-j|$ 
together with $2m-i<1$ to conclude that $i+j>2m$, and then that $i-m$ and
$j-m$ are both nonnegative, so $j \geq i$. Thus, 
$f^{j,1}(\sW_j) \subseteq f^{i,1}(\sE_i)$, so we conclude from the direct
sum decomposition of $\sE_1$ that the left side of \eqref{eq:subset} is
equal to
$(f^{i,2m-1} (\sE_i))^{\perp} \cap 
\left(\bigoplus_{j' \neq j} f_{j',1} \sW_{j'}\right)$,  
which yields the desired containment.

The cases that $2m-i >n$ and that $|m-j| \leq |m-i|$ proceed in the same 
fashion, so we conclude the lemma.
\end{proof}

\begin{rem} One might wonder whether in the definition of a linked
symplectic form, using the notation from the proof of Proposition
\ref{prop:form-structure}, it would not be enough to ask that the induced 
alternating form on $\bigoplus_i \sF_i$ be symplectic. While this condition 
might seem natural, it is not visibly intrinsic, nor does it arise naturally
from the context of limit linear series. We will see in Example
\ref{ex:symplectic} below that it is not enough to guarantee the behavior 
we want.
\end{rem}

\section{Linked symplectic Grassmannians}

We can now proceed to define linked alternating Grassmannians and 
linked symplectic Grassmannians, and we easily conclude our main 
results. 

\begin{defn} Given $\sE_{\bullet}$ $s$-linked with a linked alternating
form $\left<,\right>_{\bullet}$, the {\bf linked alternating Grassmannian}
$LAG(r,\sE_{\bullet},\left<,\right>_{\bullet})$ is the closed subscheme
of $LG(r,\sE_{\bullet})$ parametrizing linked subbundles which are 
isotropic for (the restriction of) $\left<,\right>_{\bullet}$.
\end{defn}

\begin{proof}[Proof of Theorem \ref{thm:main-alt}] By definition,
$LAG(r,\sE_{\bullet},\left<,\right>_{\bullet})$ is precisely the isotropy
locus of the restriction of $\left<,\right>_{\bullet}$ to the universal 
subbundle on $LG(r,\sE_{\bullet})$. Since the statement
is local, we may restrict to the smooth locus of $LG(r,\sE_{\bullet})$,
which according to Theorem \ref{thm:lg} is precisely the locus of exact
points. On this locus, the universal subbundle is $s$-linked, and
we conclude the desired statement from Corollary \ref{cor:iso-codim}.
\end{proof}

\begin{defn} Given $\sE_{\bullet}$ $s$-linked with a linked symplectic
form $\left<,\right>_{\bullet}$, the {\bf linked symplectic Grassmannian}
$LSG(r,\sE_{\bullet},\left<,\right>_{\bullet})$ is the closed subscheme
of $LG(r,\sE_{\bullet})$ parametrizing linked subbundles which are 
isotropic for (the restriction of) $\left<,\right>_{\bullet}$.
\end{defn}

\begin{proof}[Proof of Theorem \ref{thm:main-symp}] Once again,
$LSG(r,\sE_{\bullet},\left<,\right>_{\bullet})$ is precisely the isotropy
locus of the restriction of $\left<,\right>_{\bullet}$ to the universal 
subbundle on $LG(r,\sE_{\bullet})$, which we recall is the pullback of
the zero section under the induced morphism from $LG(r,\sE_{\bullet})$
to the space of linked alternating forms on the universal subbundle.
We may again restrict to the smooth locus of $LG(r,\sE_{\bullet})$, so
that the space of linked alternating forms is by Proposition 
\ref{prop:form-structure} a vector bundle of rank $\binom{r}{2}$, 
and we may view 
$LSG(r,\sE_{\bullet},\left<,\right>_{\bullet})$ as the intersection of
two sections inside this bundle. In order to prove the theorem, it is then
enough (see for instance Lemma 4.4 of \cite{os16}) to see that the tangent
spaces to these sections intersect transversely in the fiber over any
point of $S$. We may thus assume that $S$ is a point, and thus the
$\sE_{\bullet}$ are simply vector spaces.

At a point of the zero section, the tangent space of our bundle decomposes 
canonically as a direct
sum of the tangent space of $LG(r,\sE_{\bullet})$ (which is described
by Lemma \ref{lem:exact-tangent}) and the tangent space to the moduli space 
of linked alternating forms on the corresponding fixed linked subspace.
Since the latter moduli space is a vector space, the tangent 
space is identified with the space itself. Given a tangent vector to
$LG(r,\sE_{\bullet})$ at a point, our tautological sections yields a
tangent vector in the moduli space of linked alternating forms on the
corresponding linked subspace, which we may think of as a linked alternating
form. One checks from the definitions that if the tangent vector is
represented by $(\vp_i:\sW_i \to \sE_i/\sW_i)_i$ for some choice of
$\sW_i$ as in Lemma \ref{lem:structure}, the resulting linked alternating 
form obtained from the tautological section at this point is precisely
$\left<,\right>^{\vp_{\bullet}}_{\bullet}$ as defined in Definition
\ref{defn:tangent-form}.
Since tangent vectors to the zero section always yield the 
zero linked alternating form, transversality of the tangent spaces of
the two sections follows from the surjectivity of the map
$\vp_{\bullet} \to \left<,\right>^{\vp_{\bullet}}_{\bullet}$, given to
us by Lemma \ref{lem:symplectic-transverse}. We thus conclude the theorem.
\end{proof}

We conclude with two examples. The first demonstrates the importance of 
considering pairings between different spaces in defining a linked 
alternating form, while the second justifies our definition of a linked
symplectic form.

\begin{ex}\label{ex:need-pairings}
Consider the case $d=4$, $n=3$, $r=2$, and working over a DVR with
uniformizer $s$. We suppose we have chosen bases of the ambient spaces so
that $f_1=f_2=\begin{bmatrix} 1 & 0 & 0 &0 \\ 0 & 1 & 0 & 0\\
0 & 0 & s & 0 \\ 0 & 0 & 0 & s\end{bmatrix}$, and
that $f^1=f^2=\begin{bmatrix} s & 0 & 0 &0 \\ 0 & s & 0 & 0\\
0 & 0 & 1 & 0 \\ 0 & 0 & 0 & 1\end{bmatrix}$. The resulting linked
Grassmannian has relative dimension $4$, and we want a symplectic linked
Grassmannian to have relative dimension $3$. 

First suppose we only consider alternating forms on each individual
space, compatible with the $f_i$ and $f^i$. Then if we consider linked
alternating forms of index $2$, we could
set $\left<,\right>_{2,2} = \begin{bmatrix} 0 & 1 & 0 & 0 \\ -1 & 0 & 0 & 0 \\
0 & 0 & 0 & 1 \\ 0 & 0 & -1 & 0\end{bmatrix}$ for maximum nondegeneracy.
Compatibility with $f_1$ and $f^2$ then forces 
$$\left<,\right>_{1,1} = \begin{bmatrix} 0 & 1 & 0 & 0 \\ -1 & 0 & 0 & 0 \\
0 & 0 & 0 & s^2 \\ 0 & 0 & -s^2 & 0\end{bmatrix} \text{ and }
\left<,\right>_{3,3} = \begin{bmatrix} 0 & s^2 & 0 & 0 \\ -s^2 & 0 & 0 & 0 \\
0 & 0 & 0 & 1 \\ 0 & 0 & -1 & 0\end{bmatrix}.$$
 
Now, an open subset of linked Grassmannian can be written in the form
\begin{align*} \sF_1 & = \spn((1,a_1,0,a_2),(0,s^2 b_1, 1, b_2)),\\
\sF_2 & = \spn((1,a_1,0,s a_2),(0,s b_1, 1, b_2)),\\
\sF_3 & = \spn((1,a_1,0,s^2 a_2),(0, b_1, 1, b_2)).\end{align*}
Working over the entire DVR, we see that the condition that these subspaces
are isotropic for 
$\left<,\right>_{1,1}$, $\left<,\right>_{2,2}$, and $\left<,\right>_{3,3}$
is simply that $sb_1+sa_2=0$. Over the generic point, we get $b_1+a_2=0$,
imposing the desired additional condition. However, at $s=0$ we see that
the subspaces are automatically isotropic, so we get a full $4$-dimensional
component of the linked Grassmannian.

In order to obtain the desired relative dimension, we must impose the 
condition $b_1+a_2=0$ even on the closed fiber. We
thus see that it is necessary to consider also the pairing 
$\left<,\right>_{1,3}$, which according to our compatibility conditions will
be given by
$\left<,\right>_{1,3} = \begin{bmatrix} 0 & 1 & 0 & 0 \\ -1 & 0 & 0 & 0 \\
0 & 0 & 0 & 1 \\ 0 & 0 & -1 & 0\end{bmatrix}$. The requirement that
$\sF_1$ be orthogonal to $\sF_3$ under $\left<,\right>_{1,3}$ then
yields the desired condition $b_1+a_2=0$.
\end{ex}

\begin{ex}\label{ex:symplectic}
Consider the case $d=4$, $n=2$, $r=2$, working over a base field. We
set $s=0$, and for a given basis, consider maps of the form
$f_1=\begin{bmatrix} 1 & 0 & 0 &0 \\ 0 & 1 & 0 & 0\\
0 & 0 & 0 & 0 \\ 0 & 0 & 0 & 0\end{bmatrix}$, and
$f^1=\begin{bmatrix} 0 & 0 & 0 &0 \\ 0 & 0 & 0 & 0\\
0 & 0 & 1 & 0 \\ 0 & 0 & 0 & 1\end{bmatrix}$. As before, we obtain a
linked Grassmannian of dimension $4$, and would want a linked symplectic
Grassmannian to have dimension $3$. 

We obtain a linked alternating form of index $2$ by setting
$\left<,\right>_{2,2} = \begin{bmatrix} 0 & 0 & 0 & 1 \\ 0 & 0 & 1 & 0 \\
0 & -1 & 0 & 0 \\ -1 & 0 & 0 & 0\end{bmatrix}$, and letting the other
pairings be induced from $\left<,\right>_{2,2}$ as determined by the 
compatibility conditions. If we set $\sF_1=\spn((1,0,0,0),(0,1,0,0))$ and
$\sF_2=\spn((0,0,1,0),(0,0,0,1))$ then one checks that this induces 
a symplectic form on $\sF_1 \oplus \sF_2$. It does not satisfy our 
conditions for a linked symplectic form, because $\left<,\right>_{1,1}$
is the zero form, so its degeneracy is strictly larger than $\ker f_1$.
Correspondingly, we see that the associated
linked alternating Grassmannian is not pure of dimension $3$; indeed, it
contains the component of the linked Grassmannian on which $V_2=f_1(V_1)$.
\end{ex}

\section{Degenerations and rank-$2$ Brill-Noether loci}

We sketch how the symplectic linked Grassmannian may be used to strengthen
the limit linear series techniques of \cite{te1} to apply also in the case
of rank $2$ and canonical determinant, at least in the case that the
degenerate curve has two components. We defer a more comprehensive
exposition until the theory has been generalized to arbitrary curves of
compact type. We do, however, state a precise theorem, in terms of
the theory of \cite{te1}.

We assume throughout that $X_0$ is a
reducible projective curve of genus $g$ obtained as $Y \cup Z$, where $Y$
and $Z$ are smooth of genus $g_Y$ and $g_Z$ respectively, and 
$Y \cap Z = \{P\}$ is an ordinary node. We first recall the definition
of limit linear series in this context from \cite{te1}.

\begin{defn}\label{defn:montserrat-lls} Given $k,r,d$ positive integers, 
a {\bf limit linear series}
of dimension $k$, rank $r$, and degree $d$ on $X_0$ consists of a tuple
$(E_Y, V_Y, E_Z, V_Z, \vp_P)$, where:
\begin{ilist}
\itm $E_Y, E_Z$ are vector bundles of rank $r$ on $Y$ and $Z$ respectively; 
\itm $V_Y$ and $V_Z$ are $k$-dimensional spaces of global sections of $E_Y$ 
and $E_Z$ respectively;
\itm $\vp_P$ is an isomorphism of the projectivizations of the fibers of 
$E_Y$ and $E_Z$ at $P$,
\end{ilist}
such that there exist
\begin{ilist}[3]
\itm an integer $a>0$;
\itm bases $s^Y_1,\dots,s^Y_k$ of $V_Y$ and $s^Z_1,\dots,s^Z_k$ of $V_Z$,
\end{ilist}
satisfying the following conditions:
\begin{alist} 
\itm $\deg E_Y + \deg E_Z = d +a$;
\itm the orders of vanishing $a^Y_i$ and $a^Z_i$ of $s^Y_i$ and $s^Z_i$
at $P$ satisfy 
$$a^Y_i+a^Z_i \geq a$$
for all $i$;
\itm $s^Y_i$ glues to $s^Z_i$ under $\vp_P$ for all;
\itm global sections of $E_Y(-aP)$ and $E_Z(-aP)$ are completely determined
by their value in the fiber at $P$.
\end{alist}
\end{defn}

We then define limit linear series of canonical determinant as follows.

\begin{defn} Given $k>0$, let $(E_Y,V_Y,E_Z,V_Z,\vp_P)$ be a limit linear 
series of rank $2$, degree $2g-2$, and dimension $k$ on $X_0$. 
We say that $(E_Y,V_Y,E_Z,V_Z,\vp_P)$ has {\bf canonical determinant}
if $\det E_Y \cong \omega_Y ((d_Y-2g_Y+2)P)$, and
$\det E_Z \cong \omega_Z ((d_Z-2g_Z+2)P)$,
where $d_Y:=\deg E_Y$ and $d_Z:=\deg E_Z$.
\end{defn}

Our theorem is then the following.

\begin{thm}\label{thm:canon-det-smoothing} 
Given $g,k$ set $\rho_{\omega}=3g-3-\binom{k+1}{2}$. Suppose that
$(E_Y,V_Y,E_Z,V_Z,\vp_P)$ is a limit linear series of canonical determinant
and dimension $k$
on $X_0$ such that the inequalities of Definition \ref{defn:montserrat-lls}
b) are all equalities. Suppose further that the space of such limit linear 
series
on $X_0$ has dimension $\rho_{\omega}$ at $(E_Y,V_Y,E_Z,V_Z,\vp_P)$. Then
a general smooth curve of genus $g$ has a vector bundle of rank $2$ and
canonical determinant with at least $k$ linearly independent global 
sections.
\end{thm}

This theorem can be sharpened in a straightforward way to include
stability conditions; for these, we refer the reader to \cite{te5}.

We now sketch the proof of the theorem. Just as in \cite{os8}, we define
higher-rank limit linear series in terms of a chain of vector bundles 
$E_i$ on $X_0$ related by twisting up and down at the nodes $P$, with
maps between them given by inclusion on $Y$ and zero on $Z$ or vice versa.
A limit linear series of dimension $k$ then requires such a chain $E_i$, 
together with $k$-dimensional spaces $V_i$ of global sections of $E_i$,
each mapping into one another under the given maps. In a smoothing family
of $X_0$, the definition is the same except that the $E_i$ are related by
twisting by $Y$ and $Z$, which are now divisors on the total space. We
prove representability as in Theorem 5.3 of \cite{os8}, with moduli stacks 
of vector bundles in place of Picard schemes. Following Proposition 6.6
of \cite{os8}, we see that the forgetful map to the extremal pairs of
$E_i$ and $V_i$ yields a generalized limit linear series in the sense of
Definition \ref{defn:montserrat-lls}, with the possible exception of the 
gluing condition in part c). However, generalizing the case of refined
Eisenbud-Harris limit series, this forgetful map gives an isomorphism above
the open locus for which the inequalities of part b) are satisfied with
equality. Finally, on this locus, the corresponding points of the ambient
linked Grassmannian are all exact.

In the special case of rank $2$ and canonical determinant, one derivation 
of the modified
expected dimension for smooth curves is by constructing the moduli space
as follows. Let $\cM_{2,\omega}(X)$ be the moduli stack of vector bundles
of rank $2$ and fixed canonical determinant on $X$; this is smooth of
dimension $3g-3$. Let $\tilde{\sE}$ be the universal bundle on
$\cM_{2,\omega}(X) \times X$, and let $D$ be a sufficently ample effective
divisor on $X$ (technically, we must cover $\cM_{2,\omega}(X)$ by a nested
increasing sequence of open quasicompact substacks, and carry out this
construction on each, letting $D$ grow). Let $D'$ be the pullback of $D$ to
$\cM_{2,\omega}(X) \times X$. Then $p_{1*} \tilde{\sE}(D')$ is a vector 
bundle of rank 
$$\deg \tilde{\sE}+\rk \tilde{\sE} \deg D+ \rk \tilde{\sE}(1-g)=2g-2+2 \deg D
+2-2g= 2 \deg D.$$
Let $G:=G(k,p_{1*} \tilde{\sE}(D'))$ be the
relative Grassmannian on $\cM_{2,\omega}(X)$; our moduli space is cut out
by the closed condition of subspaces lying in $p_{1*}\tilde{\sE}$. We 
express this
condition in terms of the bundle $p_{1*}(\tilde{\sE}(D')/\tilde{\sE}(-D'))$,
which has rank $4 \deg D$. We see that because $D$ was chosen to be large,
$p_{1*} \tilde{\sE}(D')$ is naturally a subbundle, as is 
$p_{1*}(\tilde{\sE}/\tilde{\sE}(-D'))$, which also has rank $2 \deg D$. 
Then the inclusion of the universal subbundle on $G$, together with the 
pullback from $\cM_{2,\omega}(X)$ of 
$p_{1*}(\tilde{\sE}/\tilde{\sE}(-D'))$, induces a morphism
$$G \to G(k,p_{1*}(\tilde{\sE}(D')/\tilde{\sE}(-D')) \times_{\cM_{2,\omega}(X)}
G(2 \deg D,p_{1*}(\tilde{\sE}(D')/\tilde{\sE}(-D')),$$
and our desired moduli space is precisely the preimage in $G$ of the
incidence correspondence in the product. 

We now make use of the canonical determinant hypothesis to observe that
by choosing local representatives, using the isomorphism
${\bigwedge}^2 \tilde{\sE} \cong p_{2}^* \omega_X$, and summing residues
over points of $D$, we obtain a symplectic form on 
$p_{1*}(\tilde{\sE}(D')/\tilde{\sE}(-D'))$, Moreover, both 
$p_{1*} \tilde{\sE}(D')$ and $p_{1*}(\tilde{\sE}/\tilde{\sE}(-D'))$ are
isotropic for this form, with the former following from the residue theorem,
and the latter from the lack of poles. Thus, our induced map in fact
has its image in a product of symplectic Grassmannians, and the incidence
correspondence has smaller codimension, so we obtain the modified dimension
bound for our moduli space cut out in $G$.

Moving to the limit case, we need to see that the canonical determinant
hypothesis gives us (at least locally on $\cM_{2,\omega}(X))$ a linked 
symplectic form on the chain
$p_{1*}(\tilde{\sE}_i(D')/\tilde{\sE}_i(-D'))$, allowing us to extend
the above construction. We assume we have a family $X/B$ with smooth
generic fiber, and $X_0$ as above. Working locally, we may assume that each
$\tilde{\sE}_i$ has a single fixed multidegree on reducible fibers.
First suppose that for some $m$, the determinant of $\tilde{\sE}_m$ 
is isomorphic to $\omega_{X/B}$, the relative dualizing sheaf. In this
case, we construct a linked symplectic form of index $m$ by making
use of this isomorphism and our given maps between the $\tilde{\sE}_i$ to 
induce maps 
$$\tilde{\sE}_i \otimes \tilde{\sE}_j \to \omega_{X/B}$$ 
for all $i,j$, and using a fixed choice of isomorphism 
$\sO_X(Y+Z) \cong \sO_X$ to ``factor out'' any vanishing along $X_0$.
This factoring out will give us nondegeneracy whenever $i+j=2m$, and we
see that we get a linked symplectic form. On the other hand, if
no $\tilde{\sE}_i$ has determinant $\omega_{X/B}$, then for some $i$ we
have $\det \tilde{\sE}_i \cong \omega_{X/B}(Y)$. Then we have an induced
surjection
$$\tilde{\sE}_{i-1} \otimes \tilde{\sE}_i \twoheadrightarrow \omega_{X/B},$$
and if we set $m=i-\frac{1}{2},$ we get an induced linked symplectic form
as above.

Now we can put together the canonical determinant construction and the
limit linear series construction by replacing the symplectic Grassmannians
in the above argument by linked symplectic Grassmannians, and checking that
in this case we still define the same functor as before.  However, because
the linked symplectic Grassmannian is by Theorem \ref{thm:main-symp}
smooth of dimension equal to the usual symplectic Grassmannian (at least,
on the open locus of exact points), the dimension count in the limit linear
series case goes through exactly as in the case of smooth curves, and
we obtain the desired lower bound on dimension. Because the construction
goes through for smoothing families, we obtain the smoothing statement
Theorem \ref{thm:canon-det-smoothing} just as in the original Eisenbud-Harris
theory.

Note that the fact that moduli spaces of vector bundles are no longer proper
is irrelevant to us, as we are concerned only with smoothing statements.

\bibliographystyle{hamsplain}
\bibliography{hgen}

\end{document}

%% file: headers.tex
\usepackage{amssymb,euscript,amsmath, mathrsfs, amscd}
\usepackage[dvips]{graphicx}
\usepackage[dvips]{color}

\newcounter{ENUM}
\newcommand{\itm}{\item}
\newenvironment{ilist}[1][0]{\renewcommand{\theENUM}{\roman{ENUM}}\renewcommand{\itm}{\addtocounter{ENUM}{1}\item[\theENUM)]}\begin{itemize}\setcounter{ENUM}{#1}}{\end{itemize}}
\newenvironment{Ilist}{\renewcommand{\theENUM}{\Roman{ENUM}}\renewcommand{\itm}{\addtocounter{ENUM}{1}\item[(\theENUM)]}\begin{itemize}\setcounter{ENUM}{0}}{\end{itemize}}
\newenvironment{alist}[1][0]{\renewcommand{\theENUM}{\alph{ENUM}}\renewcommand{\itm}{\addtocounter{ENUM}{1}\item[\theENUM)]}\begin{itemize}\setcounter{ENUM}{#1}}{\end{itemize}}

\newcommand{\margh}[1]{}

\input xy
\xyoption{all}
\CompileMatrices

\def\ZZ{{\mathbb Z}}

\def\cM{{\mathcal M}}

\def\sE{{\mathscr E}}
\def\sF{{\mathscr F}}

\def\sO{{\mathscr O}}

\def\sW{{\mathscr W}}

\def\vp{\varphi}

\def\Hom{\operatorname{Hom}}

\def\Spec{\operatorname{Spec}}

\def\id{\operatorname{id}}

\def\rk{\operatorname{rk}}
\def\im{\operatorname{im}}

\def\spn{\operatorname{span}}
\def\sw{\operatorname{sw}}

\newtheorem{thm}{Theorem}[section]
\newtheorem{prop}[thm]{Proposition}
\newtheorem{lem}[thm]{Lemma}
\newtheorem{cor}[thm]{Corollary}

\theoremstyle{definition}
\newtheorem{defn}[thm]{Definition}

\newtheorem{ex}[thm]{Example}

\theoremstyle{remark}
\newtheorem{notn}[thm]{Notation}
\newtheorem{rem}[thm]{Remark}

\numberwithin{equation}{section}
\numberwithin{figure}{section}